\documentclass[11pt]{amsart}
\usepackage{amssymb,amsmath,amsthm,amscd, amssymb,enumerate, verbatim}
\usepackage{color}

\usepackage[all]{xy}
\numberwithin{equation}{section}
\newtheoremstyle{fancy}{10pt}{10pt}{\itshape}{12pt}{\textsc\bgroup}{.\egroup}{8pt}{
}
\newtheoremstyle{fancy2}{10pt}{10pt}{}{12pt}{\itshape}{.}{8pt}{ }

\theoremstyle{fancy}

\newtheorem{lemma}[equation]{Lemma}
\newtheorem{proposition}[equation]{Proposition}
\newtheorem{theorem}[equation]{Theorem}

\newtheorem{main}{Theorem}

\newtheorem{cor}[main]{Corollary}
\renewcommand{\themain}{\Alph{main}}
\theoremstyle{fancy2}

\newtheorem{remark}[equation]{Remark}
\newtheorem*{rem}{Remark}

\newcommand{\cref}[1]{Corollary~\ref{#1}}

\newcommand{\lref}[1]{Lemma~\ref{#1}}
\newcommand{\pref}[1]{Proposition~\ref{#1}}

\newcommand{\no}{\noindent}

\newcommand{\C}{{\mathbb{C}}}
\newcommand{\e}{\epsilon}
\newcommand{\g}{\gamma}
\newcommand{\R}{{\mathbb{R}}}
\newcommand{\Q}{{\mathbb{Q}}}
\newcommand{\CP}{{\mathbb{CP}}}
\newcommand{\HP}{{\mathbb{HP}}}
\newcommand{\Z}{{\mathbb{Z}}}
\renewcommand{\H}{{\mathbb{H}}}

\newcommand{\diag}{\ensuremath{\operatorname{diag}}}

\newcommand{\rank}{\ensuremath{\operatorname{rank}}}
\newcommand{\coker}{\ensuremath{\operatorname{coker}}}

\def\con#1=#2(#3){#1 \equiv #2 \bmod{#3}}

\begin{document}

\title{Biquotients with singly generated rational cohomology}

\author{Vitali Kapovitch}
\address{University of California\\
      Santa Barbara , CA 93110}
\email{vitali@math.ucsb.edu}
\author{Wolfgang Ziller}
\address{University of Pennsylvania\\
       Philadelphia, PA 19104}
\email{wziller@math.upenn.edu}
\thanks{Both authors were supported by grants from the
National Science Foundation. The second author was also supported
by the Francis J. Carey Term Chair.}

\maketitle

\begin{abstract}\noindent
We classify all biquotients whose rational cohomology rings are
generated by one element. As a consequence we show that the
Gromoll-Meyer $7$-sphere is the only exotic sphere which can be
written as a biquotient.
\end{abstract}

\bigskip

Let $G$ be a compact Lie group and $H\subset G\times G$ be a
compact subgroup. Then $H$ acts on $G$ on the left by the formula
$(h_1,h_2)g=h_1gh_2^{-1}$. If this action happens to be free the
orbit space is a manifold which is called a {\it biquotient} of
$G$ by $H$ and denoted by $G/\!/H$. In the special case when $H$
has the form $K_{1}\times K_{2}$ where $K_{1}\subset G\times 1$
and $K_{2}\subset 1\times G$ we will often write $K_{1}\backslash
G/K_{2}$ instead of $G//(K_{1}\times K_{2})$.

Biquotients are natural generalizations of homogeneous spaces, and
like homogeneous spaces, have  metrics with nonnegative sectional
curvature induced by  biinvariant metrics on $G$. The concept of a
biqoutient was first introduced by Gromoll-Meyer in \cite{GM},
where they showed that one of these biquotients, $Sp(2)/\!/Sp(1)$,
is an exotic 7-sphere, which produced the first example of an
exotic sphere with nonnegative curvature. Biquotients were later
on examined more systematically in \cite{Es},\cite{Bo} in the
context of a search for new manifolds with positive sectional
curvature. In fact, all known examples of manifolds admitting
metrics of positive sectional curvature are given by biquotients.

  Some
attempts were made to find other exotic spheres which could be
written as biquotients but they proved unsuccessful. We will show
in this paper that any such  attempt must indeed fail. More
generally we classify the biquotients which are rationally spheres
and projective spaces, extending a well known classification in
the homogeneous case \cite[p.195-196]{Be}:

\begin{main}
Let $M=G/\!/H$ be a compact, simply connected biquotient whose
rational cohomology ring is generated by one element. Then $M$ is
either diffeomorphic to a compact rank one symmetric space,  or it
is diffeomorphic to one of the eleven homogeneous spaces or  four
biquotients in Table \ref{mainT}.
\end{main}

\renewcommand{\thetable}{\themain}
\renewcommand{\arraystretch}{1.4}
\stepcounter{main}
\begin{table}[!ht]
         \begin{center}
             \begin{tabular}{|c|c|c|}
\hline $G/\!/H$   &\text{\rm range of $n$ }         &\text{\rm
rational type} \\
\hline\hline
$SO(2n+1)/(SO(2n-1)\times SO(2))$&$n\ge 2$&$\CP^{2n-1}$\\
\hline
  $SO(2n+1)/SO(2n-1)$&$n\ge 2$&$\mathbb{S}^{4n-1}$\\
\hline
  $SU(3)/SO(3)_2$&&$\mathbb{S}^{5}$\\
\hline
$Sp(2)/Sp(1)_{10}$&&$\mathbb{S}^{7}$\\
\hline
$G_2/SU(2)_3$&&$\mathbb{S}^{11}$\\
\hline
$G_2/SO(3)_4$&&$\mathbb{S}^{11}$\\
\hline
$G_2/SO(3)_{28}$&&$\mathbb{S}^{11}$\\
\hline
$G_2/U(2)_3$&&$\CP^{5}$\\
\hline
$G_2/SO(4)$&&$\HP^2$\\
\hline \hline
$\triangle SO(2)\backslash SO(2n+1)/SO(2n-1)$&$n\ge 2$&$\CP^{2n-1}$\\
\hline
$\triangle SU(2)\backslash SO(4n+1)/SO(4n-1)$&$n\ge 2$&$\HP^{2n-1}$\\
\hline
$Sp(2)/\!/Sp(1)$&&$\mathbb{S}^{7}$\\
\hline
$G_2/\!/SU(2)$&&$\mathbb{S}^{11}$\\
\hline
             \end{tabular}
         \end{center}
         \vspace{0.1cm}
         \caption{Rational Spheres and Projective Spaces}\label{mainT}
\end{table}
\renewcommand{\thetable}{\theequation}

Some comments may be helpful, in order to understand the examples
in this Table. The subscript for the 3 dimensional subgroups
denotes the index of the subgroup, where a simple subgroup $H$ in
a simple Lie group $G$ has index $k$ if the induced map
$\pi_3(H)\simeq\Z\to \pi_3(G)\simeq\Z$ is multiplication by $\pm
k$, which in particular means that $\pi_3(G/H)\simeq\Z_k$. Notice
that $Sp(1)_{10}$ is the unique maximal 3 dimensional subgroup in
$Sp(2)$, such that $Sp(2)/Sp(1)_{10}$ is the normal homogeneous
Berger space with positive curvature (in fact the only entry in
Table \ref{mainT} which is known to admit a metric with positive
sectional curvature).

The  subgroups in $G_2$ can be described as follows:
In $G_2$ one has the  maximal equal rank subgroups $SO(4)$ and $SU(3)$. The
subgroup $SO(4)$ contains two normal $SU(2)$'s. One of them has
index one in $G_2$ and is also contained in $SU(3)\subset G_2$.
The quotient $G_2/SU(2)_1$ is diffeomorphic to $SO(7)/SO(5)$. The
other $SU(2)\subset SO(4)$ has index 3 in $G_2$. Each $SU(2)$ can be enlarged to a $U(2)$
 $\subset SO(4)\subset G_2$. Furthermore, $G_2/U(2)_1$ is diffeomorphic to $SO(7)/SO(5)SO(2)$. One also has a
subgroup $SO(3)\subset SO(4)$ which has index 4 in $G_2$, and is
also contained in $SU(3)$. Finally there exists a maximal $SO(3)$
in $G_2$ which has index 28.

The biquotient $G_2/\!/SU(2)$ is obtained by letting $SU(2)$ act
via the index three $SU(2)$ on the left, and the index four
$SO(3)$ on the right. The Gromoll-Meyer sphere $Sp(2)/\!/Sp(1)$ is
obtained by letting $Sp(1)$ act via $\diag (q,q)$ on the left, and
$\diag (q,1)$ on the right. In the two even dimensional
biquotients, the subgroup on the left is embedded as
$\diag(1,A,\cdots ,A)$ where $A$ lies either in $SO(2)$ or in
$SU(2)$. Of these four biquotients, all but the Gromoll-Meyer
sphere were first discovered by Eschenburg in \cite{Es}, except
that he did not discuss their topological properties.

  By
computing the cohomology rings and the Pontrjagin classes, we will
show that none of these spaces are homeomorphic to each other,
except that $Sp(2)/\!/Sp(1)$ is homeomorphic to $S^7$ \cite{GM},
and the rational 11 sphere $G_2/\!/SU(2)$ is homeomorphic to
$SO(7)/SO(5)\simeq T^1S^6 \simeq G_2/SU(2)$. At the moment we are
unable to decide if the last two spaces are diffeomorphic or not,
but we at least show they can only differ from each other by a
connected sum with one of the 992 homotopy $11$-spheres.

In particular we obtain the following

\begin{cor}
The only biquotient which can be an exotic sphere is diffeomorphic
to  the Gromoll-Meyer sphere $Sp(2)/\!/Sp(1)$.
\end{cor}

Some of the other spaces given by Theorem A also have interesting
relationships. We will show that the Grassmannian
$SO(2n+1)/(SO(2n-1)\times SO(2))$ and the biquotient $\triangle
SO(2)\backslash SO(2n+1)/SO(2n-1)$ have the same cohomology rings,
and they also have the same integral cohomology groups as
$\CP^{2n-1}$, but they can be distinguished by their Pontrjagin
classes. Similarly, $\triangle SU(2)\backslash SO(4n+1)/SO(4n-1)$
has the same integral cohomology groups as $\HP^{2n-1}$, but a
different ring structure.

After a first version of this paper was finished, the preprint
\cite{T} by B.~Totaro came to our attention, where the author
independently classifies all biquotients which are rational
homology spheres. In that paper Totaro also determines exactly
which Cheeger manifolds (i.e connected sums of two compact rank
one symmetric spaces) can be written as biquotients and proves
some other interesting results about rational structure of
biquotients.

  At the same time B.~Wilking pointed out to us
that the reduction to the simple case, which in our first version
was more complicated, can be easily achieved by Lemma 1.4, which
we then noticed is also Lemma 4.5 in \cite{T}. We will use this
simplified version of the proof. We would finally like to thank
B.~ Wilking for several further useful comments.    We also thank J. DeVito for pointing out to us that
in Table B the homogeneous space $G_2/U(2)_3$ was missing in the published version of this paper.

\section {Reduction to the case of a simple $G$}\label{reduction}

Throughout this section  all cohomology have rational coefficients
and all homotopy groups are tensored with $\mathbb Q$. Also for
the purposes of the proof we will use the following equivalent but
formally stronger definition of a biquotient, see \cite{Es}:

Let $H\overset{\rho}{\to} G\times G$ be a homomorphism and let
$\Delta Z_G$ be the diagonal embedding of the center of $G$ into
$G\times G$. Let $Z=\rho^{-1}(\Delta Z_G)$. It is clear that $Z$
lies in the kernel of the usually defined biquotient action of $H$
on $G$. Suppose the biquotient action of $H/Z$ on $G$ is free.
Then the quotient space is a manifold diffeomorphic to
$(G/\rho(Z))/\!/(H/Z)$ which we will still denote by $G/\!/H$.

In this section we will also not distinguish between a simple
group and its various covers, using e.g. the same notation for
$SO(n)$ and $Spin(n)$.

We first need to recall some well known facts about rational
cohomology of Lie groups.

A Lie group of rank $n$ is rationally homotopy equivalent to a
product of a finitely many odd dimensional spheres
$S_1^{2k_1+1}\times...\times S_n^{2k_n+1}$ . The dimensions of the
spheres corresponding to various simple groups are listed in
Table~\ref{spheredim}\\

\stepcounter{equation}
\begin{table}[!ht]
         \begin{center}
             \begin{tabular}{|c|c|}
\hline $G$            &$\dim S_i$\\ \hline\hline $SO(2n-1)$
&$3,7,...,4n-5$\\ \hline $SO(2n)$  &$3,7,...,4n-5,2n-1$\\ \hline
$SU(n)$  &$3,5,...,2n-1$\\ \hline $Sp(n)$ &$3,7,...,4n-1$\\
     \hline
     $G_2$  &$3,11$\\
     \hline
     $F_4$  &$3,11,15,23$\\
     \hline
     $E_6$  &$3,9,11,15,17,23$\\
     \hline
     $E_7$  &$3,11,15,19,23,27,35$\\
     \hline
     $E_8$  &$3,15,23,27,35,39,47,59$\\
     \hline
             \end{tabular}
         \end{center}
         \vspace{0.1cm}
         \caption{Dimensions of Spheres}\label{spheredim}
\end{table}

\begin{proposition}\label{p:reduction}
Suppose $M=G/\!/H$ is a biquotient such that $M$ is simply
connected and the cohomology algebra $H^*(M,\mathbb Q)$ is
generated by one element. Then there exists a biquotient
$G'/\!/H'$ such that $G'$ is simple and $M$ is diffeomorphic to
$G'/\!/H'$.
\end{proposition}

\begin{proof}

The key in the proof of \pref{p:reduction} is the following
elementary Lemma the proof of which is left to the reader.

\begin{lemma}\label{l:reduction} Let $G$ be a compact Lie group
acting differentiably on manifolds $X$ and $Y$. Suppose that the
action of $G$ on $X$ is transitive and the diagonal action of $G$
on $X\times Y$ is free. Then for any $x\in X$ the action of
isotropy group $G_x$ on $Y$ is free and the quotient spaces
$(X\times Y)/G$ and $Y/G_x$ are canonically diffeomorphic.
Moreover, if the action of $G$ on $X\times Y$ is a biquotient
action then the action of $G_x$ on $Y$ is again a biquotient
action.
\end{lemma}

First notice that by passing to a finite cover we can assume that
both $G$  and $H$ are  products of compact simple or abelian
groups. Indeed, let $\pi :G' \to G$ be a finite cover of $G$ that
splits as a product of simple or abelian groups
$G'=G_1\times...\times G_n$. Let $\hat{H}=\pi^{-1}(H)$. Then
$G'/\!/\hat{H}\simeq G/\!/H$ and since $M$ is simply connected,
$\hat{H}$ is connected. Let $H'\to \hat{H}$ be a finite cover such
that $H'$ splits as $H'=H_1\times..\times H_m$. Then
$G'/\!/H'\simeq G'/\!/\hat{H}\simeq G/\!/H$. From now on we can
assume that $G$ and $H$ already have the product forms
$G=G_1\times...\times G_n, H=H_1\times..\times H_m$. Furthermore,
since $M$ is simply connected, by~\lref{l:reduction} we can assume
that $G$ has no abelian factors.

Next let us describe the rational homotopy type of $M$. We are
given that $H^*(M,\mathbb Q)$ is generated by one element $a$,
which easily implies that $M$ is formal. Indeed, the naturally
defined map $(H^*(M,\R),0)\to (\Omega^*(M),d)$ is clearly a DGA
quasi-isomorphism. Thus $M$ is formal over $\R$ and hence, by the
field extension theorem ~\cite[page 156]{FHT}, it is formal over
$\Q$.

     If
$\deg a$ is odd, it is obvious that $\dim M=\deg a$ and $M$ is
rationally equivalent to $S^{\deg a}$. If $\deg a =2k$ is even,
$H^*(M)=\mathbb Q[a]/a^{m+1}$ where $m=\dim M/\deg a$. By
formality, the minimal model of $M$ is then the same as the
minimal model of $(\mathbb Q[a]/a^{m+1},0)$ and is equal to
$(\mathbb Q[x, y],d)$ where $\deg x= \deg a,\deg y=(m+1)\deg a -1,
dx=0, dy=x^{m+1}$. In particular, $M$ has exactly two nontrivial
rational homotopy groups $\pi_{\deg a}(M)\backsimeq \pi_{\deg
y}(M)\backsimeq \mathbb Q$. In either case $M$ has exactly one
nontrivial odd homotopy group.

Let $i=(i_1,\ldots,i_n):H\to G^2= G_1^2\times\ldots G_n^2$ be the
fiber inclusion.

By looking at the long exact homotopy sequence of the fibration
$H\overset{i}{\to} G\to M$ we see that the induced map
$i_\ast:\pi_\ast (H)\to \pi_\ast (G)$ satisfies $\dim \coker
(i_\ast)=1$ and $\dim \ker (i_\ast)=1 (=0)$ if $\dim M$ is even
(odd). Since $\rank(G)=\dim \pi_\ast (G)$ and $\rank(H)=\dim
\pi_\ast (H)$ this implies that $\rank (G)=\rank(H)$ if $\dim M$
is even and $\rank (G)=\rank(H)+1$ if $\dim M$ is odd.

By the above, all but one of the coordinate projections $i_k:H\to
G_k$ are onto on $\pi_\ast$.
\begin{lemma}\label{l:inj}
Let $f:H\to G$ be a continuous map between compact connected Lie
groups. Suppose the induced map $f_\ast:
\pi_\ast(H)\to\pi_\ast(G)$ is onto. Then $f$ is onto.
\end{lemma}
\begin{proof}
We are going to show that the induced map $f^\ast:H^*(G)\to
H^*(H)$ is injective. First observe that since both $H$ and $G$
are rationally products of odd-dimensional spheres their
cohomology algebras are  free exterior algebras on a finite number
of odd-dimensional generators.  Thus for both $H$ and $G$ the
vector spaces  spanned by those generators  (denoted by $V_H$ and
$V_G$ respectively) can be naturally identified with quotients of
$H^*$ by decomposable elements $H^{*+}/(H^{*+}\cdot H^{*+})$. The
assumptions of the Lemma imply that the induced map
$f^*:H^{*+}(G)/(H^{*+}(G)\cdot H^{*+}(G))\to
H^{*+}(H)/(H^{*+}(H)\cdot H^{*+}(H))$ is injective. Since
$H^*(H)\simeq \Lambda V_H$ this implies that the map
$f^\ast:H^*(G)\to H^*(H)$ is injective. In particular, the image
of the fundamental cohomology class $[G]$ is nonzero and hence $f$
is onto.
\end{proof}
  By~\lref{l:inj} for all but one factor $G_i$ the action of $H$ on
  $G_i$ is transitive. Therefore by \lref{l:reduction}
   we can reduce the number of simple factors of $G$ to one.
This concludes the proof of
  ~\pref{p:reduction}.

\end{proof}

\section{Case of a simple $G$ and Proof of Theorem A}

We are now ready to proceed with the proof of Theorem A in the
Introduction. We can assume that $H^*(M)=\mathbb Q[a]/a^{m+1}$,
and that $M=G/\!/H$ with $G$ simple.

Let the embedding $H\subset G\times G$ be given by $(j^-,j^+)$
where $j^-$ and $j^+$ are two homomorphisms. If one of these is
trivial, we are in the situation of a homogeneous space where we
can use the classification in \cite[p.195-196]{Be} or \cite{On} to
obtain the first half of Table \ref{mainT}.

We will now distinguish between the case $\dim M$ odd and $\dim M$
even and use results from the proof of~\pref{p:reduction} in each
case, as well as the fact that for any simple Lie group $\rank
\pi_3=1$.

If $\dim M$ is odd and hence $H^\ast(M)=H^\ast(S^{2n+1})$, we have
$H^*(G)\cong H^*(H\times S^{2n+1})$ as rings and $\rank G=\rank
H+1$. If $\dim M =3$, $H$ must be trivial and hence $G/\!/H$ is
homogeneous. If $\dim M
>3$, then $G$ and $H$ are both simple, and hence $j^-$ and $j^+$
are either homomorphisms with finite kernels,
or trivial. Now one can easily produce a list of all simple pairs
$H\subset G$ such that $H^*(G)\cong H^*(H\times S^{2n+1})$, using
Table \ref{spheredim} and elementary representation theory. The
result is summarized in Table \ref{t:homog}. Notice that this
happens to agree with the list of homogeneous spaces $G/H$ which
are odd dimensional rational homology spheres (see
\cite[p.195-196]{Be} and \cite{On}), although this is not a priori
clear.

\stepcounter{equation}
\begin{table}[!ht]
         \begin{center}
             \begin{tabular}{|c|c|c|c|}
\hline  $G$  &$H$&range of $n$& number of reps\\ \hline\hline

$SO(2n)$ &$SO(2n-1)$&$n\ge 3, n\ne 4$& 1\\ \hline

$SU(n)$ &$SU(n-1)$&$n\ge 4$&1\\ \hline

$Sp(n)$&$Sp(n-1)$&$n\ge 3$&1\\ \hline

$SO(2n+1)$&$SO(2n-1)$&$n\ge 3$&1\\ \hline

$Spin(7)$  & $G_2$&&1\\ \hline

$Spin(8)$ &$Spin(7)$&& 3\\ \hline

$Spin(9)$ &$Spin(7)$&& 2\\ \hline

$SU(3)$   & $SU(2)$&&2\\ \hline

$Sp(2)$ & $Sp(1)$&&3\\ \hline

$G_2$   & $SU(2)$&&4\\ \hline

\end{tabular}
         \end{center}
         \vspace{0.1cm}
         \caption{Rational odd dimensional homology spheres}\label{t:homog}
\end{table}

In the first 5 cases, the embedding of $H$ in $G$ is unique up to
conjugacy and hence these cases only give rise to homogeneous
biquotients. In the remaining cases there exist at least two
embeddings of $H$ and hence the possibility of a biquotient.

The three representations of $Spin(7)$ in $Spin(8)$, as well as
the two representations of $Spin(7)$ in $Spin(9)$ intersect in
$G_2$ and hence this case cannot give rise to a biquotient. The
group $SU(3)$ has the index 1 subgroup $SU(2)$ and the index 2
subgroup $SO(3)$ which intersect in a circle and hence cannot give
rise to a biquotient. The group $Sp(2)$ has the index one subgroup
$Sp(1)\times 1\subset Sp(1)\times Sp(1)\subset Sp(2)$, the index 2
subgroup $\triangle Sp(1)\subset Sp(1)\times Sp(1)\subset Sp(2)$,
and the maximal index 28 subgroup $Sp(1)\subset Sp(2)$. It is not
hard to see that only the first two can be combined to give rise
to a biquotient, the Gromoll Meyer sphere $Sp(2)/\!/Sp(1)$. The
exceptional group $G_2$ has  4 three dimensional subgroups
described in the Introduction. The question which biquotients this
gives rise to is more complicated. However, the general situation
of $\rank G=2$ and $\rank H=1$ has been completely examined in
\cite[p.166-170]{Es} where it was shown that it gives rise to only
two biquotients. The first one is the Gromoll Meyer sphere and the
second one is $G_2/\!/SU(2)$, where one uses the index 3 and index
4 subgroups for $j^-$ and $j+$.

If $\dim M$ and hence $\deg a$ is even, we have $\rank G = \rank
H$, and since $G$ is simple, it follows that $H$ is simple if
$\deg a> 4$, $H$ has two simple factors if $\deg a =4$, and
$H=H_1\times S^1$ with $H_1$ simple if $\deg a=2$.

If $G/\!/H$ is not homogeneous, the maximal torus in $H$ must give
rise to a (two-sided) biquotient action of a torus on a simple Lie
group $G$, whose dimension is equal to the rank of $G$. These were
all classified in \cite{Es}. Such biquotient tori actions are
fairly rare, and in particular none exist for the exceptional Lie
groups. Furthermore these tori are all such that there exists a
codimension 1 torus which acts only on one side of $G$, say on the
right, and the remaining circle either acts on the left or on both
sides. Hence it follows that the image of the projection of $H$
onto the right side is a rank one group and the kernel a $\rank H
-1$ normal subgroup of $H$. Hence $H=H_1\times H_2$ with $H_2$
simple and $H_1$ is either $S^1$, $SO(3)$ or $SU(2)$. Hence
$G/H_2$ must be a homogeneous space which is  either an odd
dimensional rational homology sphere or is rationally equivalent
to $M\times H_1$. Since both $G$ and $H_2$ are simple it easily
follows that $G/H_2$ must be a rational sphere. Now we use Table
\ref{t:homog} and determine if a further rank 1 group $H_1$ can
act freely on it. If $G/H_2$ is diffeomorphic to a sphere, the
action of $G$ on this sphere is linear and hence $H_1$ can only be
one of the Hopf actions which implies that the quotient is
diffeomorphic to a projective space.

According to Table \ref{t:homog}, the only remaining cases are
$G/H_2= SO(2n+1)/SO(2n-1)$ or $SU(3)/SO(3) , Sp(2)/Sp(1)_{10}  $
and $ G_2/SU(2)$, where we have used the fact that $Sp(2)/\Delta
Sp(1)=SO(5)/SO(3)$. In each case we now have to determine if $H_1$
can act freely on it.  But Eschenburg's classification of maximal
tori that act freely immediately implies that the only
possibilities are the two entries  $\triangle SO(2)\backslash
SO(2n+1)/SO(2n-1)$ and $\triangle SU(2)\backslash
SO(4n+1)/SO(4n-1)$ in Table \ref{mainT}. Here $\triangle SO(2)$
and $\triangle SU(2)$ stand for "Hopf actions" $\diag(1,A,\cdots
,A)$ where $A$ lies either in $SO(2)$ or in $SU(2)$.  This
finishes the proof of Theorem A.

\begin{rem}
In the homogeneous case, one not only has a diffeomorphism
classification of the homogeneous spaces which are rational
homology spheres, but in each case can also determine in how many
ways the manifold can be written as a homogeneous space. In the
case of biquotients, such a classification is possible if $G$ is
simple, using \cite[Table 101]{Es}. Notice that in this Table,
there are quite a few biquotients which are diffeomorphic to $\C
P^n$ or $\H P^n$ without being homogeneous. If $G$ is not simple,
there are many possibilities, and in fact we can increase the
number of simple factors in $G$ arbitrarily by using the fact that
$G//H=\Delta G \backslash G\times G/H$ repeatedly.
\end{rem}

\section{Diffeomorphism classification}
\begin{theorem}\label{t:diff}
None of the spaces listed in Theorem A are mutually diffeomorphic
except possibly the rational $11$-spheres $G_2/\!/SU(2)$ and
$SO(7)/SO(5)$ (which can also be written as $G_2/SU(2)_1$). These
two spaces are PL-homeomorphic but may possibly differ by a
connected sum with an exotic $11$-sphere.
\end{theorem}
\begin{proof}
The homogeneous spaces can easily be differentiated from the rank
one symmetric spaces and from each other by the torsion in their
cohomology, see e.g. \cite{MZ}.

Here we will only need the integral cohomology groups of
$SO(2n+1)/SO(2n-1)=T^1S^{2n}$, which follows easily from the Gysin
sequence of the bundle $S^{2n-1}\to T^1S^{2n}\to S^{2n}$ :
\[
H^*(T^1S^{2n})=\begin{cases}\Z \text{ if }*=0,4n-1\\Z_2 \text{ if
} *=2n\\ 0 \text{ otherwise}
\end{cases}
\]

  Let us next consider the rational
$\mathbb{CP}^{2n-1}$'s $M=\triangle SO(2)\backslash
SO(2n+1)/SO(2n-1)$ and $N=SO(2n+1/SO(2n-1)\times SO(2), n>1$.

Let us first compute the integral cohomology rings of $M$ and $N$.
 From the Gysin sequence of the bundle $S^1\to T^1S^{2n}\to M$ we
compute

\[
H^*(M)=\begin{cases}\Z \text{ if }*= 2k, \text{ for
}k=0,...,2n-1\\ 0 \text{ otherwise}
\end{cases}
\]
Moreover from the same sequence we see that the Euler class of
this bundle $a_M\in H^2(M)$ is a generator of $H^2(M)$ and the
following sequences are exact

\begin{align*}& 0\to H^{2k-2}(M)\overset{\cup a_M}{\to}H^{2k}(M)\to 0
\text{ for } k=1,...,n-1,n+1,...,2n-1\\&0\to
H^{2n-2}(M)\overset{\cup a_M}{\to}H^{2n}(M)\to Z_2\to 0
\end{align*}

Hence the ring structure is determined by the fact that $a_M^n$ is
twice a generator in $H^{2n}(M)$ and thus $M$ is not homotopy
equivalent to $\C P^{2n-1}$.

The same argument works for $N$ and thus $M$ and $N$ have
isomorphic cohomology rings. To compare $M$ and $N$ to each other
we will show that they have different rational Pontrjagin classes
and thus  are not homeomorphic.

Observe that both $M$ and $N$ are quotients of $T^1(S^{2n})$ by
different free $S^1$ actions. Let us describe these actions
explicitly. We will identify $T^1(S^{2n})$ with the set\\
$\{(x,y)\in \mathbb {R}^{2n+1}\times \mathbb {R}^{2n+1}|
|x|=|y|=1, \langle x, y\rangle =0\}$. By construction, the $S^1$
action producing $M$ is the diagonal action $z(x,y)=(z(x),z(y))$
for the embedding $S^1\to SO(2n)\to SO(2n+1)$ with the first
embedding given by the Hopf action. Observe that this action
leaves the product $S^{2n}\times S^{2n}$ invariant. It is easy to
see that the normal bundle $\nu$  of $T^1(S^{2n})$ inside $\mathbb
{R}^{2n+1}\times \mathbb {R}^{2n+1}$ is trivial. Consider the
natural orthonormal trivialization $e:T^1(S^{2n})\times
\mathbb{R}^3\to \nu$ given by $e_1(x,y)=e( (x,y),(1,0,0))= (x,0),
e_2(x,y)=e( (x,y),(0,1,0))= (0,y),e_3(x,y)=e( (x,y),(0,0,1))=
\frac{1}{\sqrt{2}}(y,x)$. It is easy to see that with respect to
$e$  the action of $S^1$ on $\mathbb{R}^3$ is trivial and
therefore $\nu$ descends to a trivial bundle over $M$.

Let $p: T^1(S^{2n})\to M$ be the canonical projection. Then $T
T^1(S^{2n})\simeq p^\#(TM)\oplus T_F$ where $T_F$ is the tangent
bundle to the fiber. It is obvious that $T_F\simeq \e^1$ is a
trivial bundle over $T^1(S^{2n})$ and therefore
\[
T T^1(S^{2n})\simeq p^\#(TM)\oplus \e^1
\]

Next note that the action of $S^1$ on $\mathbb {R}^{2n+1}\times
\mathbb {R}^{2n+1}$ is equivalent to the sum of $2n$ copies of the
standard representation and a $2$-dimensional trivial
representation $\e^2$.

Combining the previous formulas we obtain the following identity
\begin{equation}\label{TM}
TM\oplus \e^4\simeq 2n\gamma_M\oplus \e^2
\end{equation}

where $\g_M$ is the rank-two bundle over $M$ associated to the
principal $S^1$ bundle $T^1(S^{2n})\to M$ and the canonical
$S^1\simeq SO(2)$ action on $\R^2$. By construction, $
e(\g_M)=a_M$, the generator of $H^2(M)\simeq \Z$.

Therefore,
\begin{equation}\label{p_M}
p_1(TM)=p_1(TM\oplus
\e^4)=p_1(2n\g_M)=2np_1(\g_M)=2ne(\g_M)^2=2na_M^2
\end{equation}

Let us now compute the first Pontrjagin class of $N$. By
definition, the $S^1$ action on $T^1(S^{2n})$ which produces $N$
is given by the following formula:

\[
e^{it}(x,y)=(\cos t x+\sin t y, -\sin t x+\cos t y)
\]
In other words this is just the geodesic flow action for the round
metric on $S^{2n}$.

As before we see that it is equivalent to the sum of $2n+1$ copies
of the standard representation and therefore it descends to the
bundle $(2n+1)\g_N$ over $N$.

On the other hand, by the same argument as before we see that
\begin{equation}\label{tan1}
(2n+1)\g\simeq TN\oplus \bar{\nu}\oplus \e^1 \end{equation}
  where $\bar{\nu}$ is
the $S^1$ quotient of the normal bundle $\nu$ to $T^1(S^{2n})$
inside $\mathbb {R}^{2n+1}\times \mathbb {R}^{2n+1}$. Let us study
$\bar{\nu}$ further. As was discussed earlier, $\nu$ is a trivial
bundle. It is easy to see that with respect to the trivialization
$e=(e_1,e_2,e_3)$ the action $\rho$ of $S^1$ on $\R^3$
corresponding to $\nu$ is given by the following matrix

\begin{equation}\label{nu}
e^{it}\longrightarrow
\begin{pmatrix}
\cos^2 t&\sin^2 t&\sqrt{2}\sin t\cos t\\ \sin^2 t&\cos^2
t&-\sqrt{2}\sin t\cos t\\ -\sqrt{2}\sin t\cos t&\sqrt{2}\sin t\cos
t&\cos^2t-\sin^2t
\end{pmatrix}
\end{equation}

 From formula~\ref{nu} we see that $\rho$ is equivalent to the sum
of a rank-one trivial representation and a representation of
weight $2$. Therefore $\nu$ descends to the bundle
$\eta\oplus\e^1$ where $\eta_N$ is a rank-two bundle over $N$ with
Euler class $2a_N$. We can now rewrite  formula~\ref{tan1} as
follows
\begin{equation}\label{tan2}
(2n+1)\g\simeq TN\oplus\eta_N\oplus\e^2
\end{equation}
Therefore
$(2n+1)a_N^2=(2n+1)p_1(\g_N)=p_1(TN)+p_1(\eta_N)=p_1(TN)+e(\eta_N)^2=p_1(TN)+4a_N^2$
and hence
\begin{equation}\label{p_N}
p_1(TN)=(2n-3)a_N^2
\end{equation}

Finally, observe that the groups $H^4(M)/\langle p_1(M)\rangle$
and $H^4(N)/\langle p_1(N)\rangle$ are  cyclic. By comparing
(\ref{p_M}) and (\ref{p_N}) we see that these groups have
different orders and therefore $M$ and $N$ are not homeomorphic by
topological invariance of rational Pontrjagin classes.

Next let us consider the rational $\HP^{2n-1}$ given by
$M=\triangle SU(2)\backslash SO(4n+1)/SO(4n-1)$. A similar
computation to the one in case of rational $\CP^n$'s shows that it
has the following cohomology
\[
H^*(M)=\begin{cases}\Z \text{ if }*= 4k, \text{ for
}k=0,...,2n-1\\ 0 \text{ otherwise}
\end{cases}
\]

\no Also, as before,  if $a_M\in H^4(M)$ is a generator of
$H^4(M)$ then the following sequences are exact

\begin{align*}& 0\to H^{4k-4}(M)\overset{\cup a_M}{\to}H^{2k}(M)\to 0
\text{ for } k=1,...,n-1,n+1,...,2n-1\\&0\to
H^{4n-4}(M)\overset{\cup a_M}{\to}H^{4n}(M)\to Z_2\to 0
\end{align*}
Therefore $a_M^n$ is twice a generator in $H^{4n}(M)$ and hence
$M$ is not homotopy equivalent and hence not diffeomorphic  to
$\HP^{2n-1}$.

The biquotient $Sp(2)/\!/Sp(1)$ is homeomorphic but not
diffeomorphic to $\mathbb{S}^7$ according to~\cite{GM}.

Let us finally discuss the rational $11$-sphere
$M^{11}=G_2/\!/SU(2)$. Recall that
\[
H^\ast(G_2)=\begin{cases}\Z \text{ if }\ast=0,3,11,14\\Z_2 \text{
if } \ast=6,9\\ 0 \text{ otherwise}
\end{cases}
\]

Since the fiber in the fibration  $SU(2)\to G_2/\!/SU(2)$ is given
by the
   composition of two maps $(j^-,j^+): SU(2) {\hookrightarrow}
G_2\times G_2$ and
   $\times :G_2\times G_2\to G_2$, where $\times(g_1,g_2)=g_1\cdot
   g_2^{-1}$, it induces the map $j^-_\ast - j^+_\ast$ in $\pi_3$.
Since $j^-$ is given by the index 3 subgroup and
   $j^+$ by the index 4 subgroup, it follows that
the fiber inclusion $SU(2)\to G_2$ is an isomorphism on $\pi_3$.
>From the long exact homotopy sequence of the fibration
$S^3=SU(2)\to G_2\to M$ we conclude that $M$ is $4$-connected.
Therefore the Euler class $e\in H^4(M)$ of this bundle is zero.
 From the Gysin sequence
\[
\to H^1(M)\overset{\cup e}{\to} H^5(M)\to H^5(G_2)\to
\]
we see that $H^5(M)=0$. Similarly, from
\[
\to  H^2(M)\overset{\cup e}{\to} H^6(M)\to H^6(G_2)\to H^3(M)\to
\]
we see that $H^6(M)\simeq H^6(G_2)\simeq \Z_2$. Thus
\[
H^\ast(M)=\begin{cases}\Z \text{ if }\ast=0,11\\Z_2 \text{ if }
\ast=6\\ 0 \text{ otherwise}
\end{cases}
\]
and by Poincare duality
\begin{equation}\label{ts6}
H_\ast(M)=\begin{cases}\Z \text{ if }\ast=0,11\\Z_2 \text{ if }
\ast=5\\ 0 \text{ otherwise}
\end{cases}
\end{equation}

We will say that two closed  manifolds are almost diffeomorphic if
they  differ by a connected sum with a homotopy sphere.
\begin{lemma} \label{almostdif}Suppose $X^{11}$ is a simply-connected
smooth manifold with homology given by~(\ref{ts6}).
  Then $X$ is almost diffeomorphic to $T^1S^6$.
\end{lemma}
\begin{proof}
The almost  diffeomorphism   classification of $k$-connected
$2k+1$ manifolds with $k\ne 3,7$ was carried out by
Wall~\cite{Wa}. By assumptions, our manifold is as above with
$k=5$.

  According to Wall,  the {\it oriented} almost
diffeomorphism class of a $4$-connected $11$-manifold is
completely determined by the following set of invariants:

\begin{enumerate}[$\bullet$]
\item $G=H_5(X)$; \item A nonsingular bilinear form  (called the
linking form) $b: G^\ast\times G^\ast\to G^\ast$where $G^*$ is the
torsion subgroup of $G$; \item A quadratic form $q: G^\ast \to
\Q/2\Z$ associated with the bilinear form $2b$; \item A
homomorphism $\alpha: G\to \pi_4(SO)$.
\end{enumerate}

In our case $G\simeq G^\ast\simeq\Z_2$. Since there exists only
one non-degenerate bilinear form on $\Z_2$, the form $b$ is
uniquely determined.

By Bott periodicity, $\pi_4(SO)=0$, and thus $\alpha =0$. The only
remaining Wall invariant which has to be determined is the
quadratic form $q$. To compute it we need to recall its
definition.

  Look at the generator $a$ in $H_5(X)=Z_2$ given by a map $a:S^5\to X$.
  By Whitney's  theorem we can assume that $a$ is an
embedding. The normal bundle to $a$ is trivial since
$\pi_4(SO(6))=0 $.

Choose a section $a_1$ of the normal bundle to $a$ such that the
normal bundle to $a _1$ in the unit tangent bundle of $a$ is
trivial. This is not automatic since $\pi_4(SO(5))=Z_2$ . The
easiest way to achieve this is to take the obvious section
corresponding to any trivialization of the normal bundle. Let $a_2
:S^5\to M$ be the normal sphere in the unit tangent bundle. The
orientation on $a_2$ is uniquely determined by the orientations on
$M$ and $a$. More explicitly, we orient the normal $D^6$ to have
intersection with $a$ equal to $+1$ and consider the induced
orientation on the normal $S^5=\partial D^6$. Let $Y=X\backslash
a(S^5)$ and let $y_1=[a_1],y_2=[a_2]$ be the homology classes in
$X$ given by $\alpha_i$.  Then it can be shown that $y_2$
generates the kernel of the map $H_5(Y) \to H_5(X)$ which is
infinite cyclic. It is clear that $2y_1$ lies in that kernel and
therefore $2y_1=\lambda y_2$. It can be shown ~\cite{Wa} that the
quotient $\lambda/2$ is well-defined mod $2\Z$ and
  we set $q(a):\overset{\text{def}}{=}\lambda/2$ mod $2$. It is
obvious that $\lambda$ can not be even so $q(a)$ can only take
values  $\pm 1/2$ mod $2$.

If we change the orientation of
  $X$ then, by construction, $y_2$
changes to $-y_2$ and hence $q$ changes to $-q$. Therefore, $X$
and $-X$  ( which stands for the same manifold with opposite
orientation) have different oriented almost diffeomorphism types
and any other oriented manifold  satisfying the assumptions of the
Lemma (e.g. $T^1S^6$) is orientably almost diffeomorphic to either
$X$ or $-X$.
\end{proof}
\begin{remark} Observe that two $11$-manifolds
satisfying~\lref{almostdif} are almost diffeomorphic iff they are
$PL$-equivalent. Indeed,  a manifold $X$
satisfying~\lref{almostdif} is homeomorphic to $T^1S^6$. In
particular it admits a $CW$ decomposition $e^0\cup e^5\cup e^6\cup
e^{11}$. Look at the universal bundle $PL/O\to B_O\to B_{PL}$. We
wish to classify different smooth structures inside a fixed $PL$
structure on $X$, i.e we have to classify the homotopy types of
all possible lifts of the classifying map $f:X\to B_{PL}$. By the
general obstruction theory, the obstruction $o_i$ to extend a
homotopy between two lifts from the $i-1$'th to the $i$'th
skeleton of $X$ lives in $H^{i}(X,\pi_i(O/PL))$.  Since $O/PL$ is
$6$-connected we have $o_i=0$ for $i=1,\ldots,6$. The CW structure
of $X$ then implies that $o_i=0$ for $i=7,\ldots,10$. Thus the
only possible nontrivial obstruction is $o_{11}$ and the PL class
of $X$ contains at most
$|H^{11}(X,\pi_{11}(O/PL))|=|\pi_{11}(O/PL)|=992$ distinct
diffeomorphism types. On the other hand, connected sums of $X$
with different homotopy spheres have different Eells-Kuiper
invariants~\cite{EK} and thus are non-diffeomorphic. Since there
are exactly   $|\pi_{11}(O/PL)|=992$ homotopy $11$-spheres, the
conclusion follows.

Thus the  oriented diffeomorphism type of $M$ is determined by its
oriented $PL$-homeomorphism type together with the Eells-Kuiper
invariant of $M$ which at the moment we are unable to compute.
\end{remark}
\end{proof}

\providecommand{\bysame}{\leavevmode\hbox
to3em{\hrulefill}\thinspace}

\end{document}